\documentclass[12pt]{amsart}
\usepackage{amssymb}
\usepackage{amsmath}
\usepackage[active]{srcltx}
\usepackage{t1enc}
\usepackage[latin2]{inputenc}
\usepackage{verbatim}
\usepackage{amsmath,amsfonts,amssymb,amsthm}
\usepackage[mathcal]{eucal}
\usepackage{enumerate}
\usepackage[centertags]{amsmath}
\usepackage{graphics}

\newtheorem{theorem}{Theorem}

\newtheorem{lemma}{Lemma}

\begin{document}
\author{G. Tephnadze}
\title[Strong summability]{A note on strong summability of two-dimensional
Walsh-Fourier series}
\address{George Tephnadze, The University of Georgia, School of Science and Technology, 77a Merab Kostava St, Tbilisi, 0128,  Georgia.}
\email{g.tephnadze@ug.edu.ge}
\thanks{The research was supported by Shota Rustaveli National Science
Foundation grant YS-18-043.}
\date{}
\maketitle

\begin{abstract}
In this paper we investigate strong summability of  the two-dimensional Walsh-Fourier series obtained in Weisz \cite{We} (see Theorem W) and prove
sharpness of this result.
\end{abstract}

\textbf{2010 Mathematics Subject Classification.} 42C10.

\textbf{Key words and phrases:} Walsh system, two-dimensional Walsh-Fourier series, Strong summability, martingale Hardy space.

Let $\mathbb{N}_+$ denote the set of positive integers, $\mathbb{N}:=\mathbb{%
N}_+ \cup \{0\}$. Denote by $Z_{2}$ the discrete cyclic group of order 2,
that is $Z_{2}=\{0,1\},$ where the group operation is the modulo 2 addition
and every subset is open. The Haar measure on $Z_{2}$ is given such that the
measure of a singleton is 1/2. Let $G$ be the complete direct product of the
countable infinite copies of the compact groups $Z_{2}.$ The elements of $G$
are of the form $x=\left( x_{0},x_{1},...,x_{k},...\right) $ with $x_{k}\in
\{0,1\}\left( k\in \mathbb{N}\right) .$ The group operation on $G$ is the
coordinate-wise addition, the measure (denote\thinspace $\,$by$\,\,\mu $)
and the topology are the product measure and topology of $Z_{2}=\{0,1\}$. The compact Abelian
group $G$ is called the Walsh group. A base for the neighborhoods of $G$ can
be given in the following way: 
\begin{eqnarray*}
I_{0}\left( x\right) &:&=G,\,\,\,I_{n}\left( x\right) :=\,I_{n}\left(
x_{0},...,x_{n-1}\right) \\
&:&=\left\{ y\in G:\,y=\left( x_{0},...,x_{n-1},y_{n},y_{n+1},...\right)
\right\} ,\,\left( x\in G,n\in \mathbb{N}\right) .
\end{eqnarray*}
These sets are called the dyadic intervals. Let $0=\left( 0:i\in \mathbb{N}%
\right) \in G$ denote the null element of $G,\,\,\,I_{n}:=I_{n}\left(
0\right) \,\left( n\in \mathbb{N}\right) .$ Set $e_{n}:=\left(
0,...,0,1,0,...\right) \in G$ the $n\,$th\thinspace coordinate of which is 1
and the rest are zeros $\left( n\in \mathbb{N}\right) .$ Let $\overline{I}%
_{n}:=G\backslash I_{n}.$

If $n\in \mathbb{N}$, then $n=\sum\limits_{i=0}^{\infty }n_{i}2^{i},$ where $%
n_{i}\in \{0,1\}\,\,\left( i\in \mathbb{N}\right) $, i.e. $n$ is expressed
in the number system of base 2. Denote $\left| n\right| :=\max \{j\in 
\mathbb{N}:n_{j}\neq 0\}$, that is, $2^{\left| n\right| }\leq n<2^{\left|
n\right| +1}.$

It is easy to show that for every odd number $n_{0}=1$ and we can write $%
n=1+\sum_{i=1}^{\left\vert n\right\vert }n_{j}2^{i}$, where $n_{j}\in
\left\{0,1\right\} ,$ $~(j\in \mathbb{N}_+)$.

For $k\in \mathbb{N}$ and $x\in G$ let as denote by 
\begin{equation*}
r_{k}\left( x\right) :=\left( -1\right) ^{x_{k}}\,\,\,\,\,\,\left( x\in
G,k\in \mathbb{N}\right)
\end{equation*}%
the $k$-th Rademacher function.

The Walsh-Paley system is defined as the sequence of Walsh-Paley functions: 
\begin{equation*}
w_{n}\left( x\right) :=\prod\limits_{k=0}^{\infty }\left( r_{k}\left(
x\right) \right) ^{n_{k}}=r_{\left\vert n\right\vert }\left( x\right) \left(
-1\right) ^{\sum\limits_{k=0}^{\left\vert n\right\vert
-1}n_{k}x_{k}}\,\,\,\,\,\,\left( x\in G,n\in \mathbb{N}_+\right) .
\end{equation*}

The Walsh-Dirichlet kernel is defined by 
\begin{equation*}
D_{n}\left( x\right) =\sum\limits_{k=0}^{n-1}w_{k}\left( x\right) .
\end{equation*}

Recall that (see $\left[ 8,\text{ p. 7}\right] $) 
\begin{equation}
D_{2^{n}}\left( x\right) =\left\{ 
\begin{array}{c}
2^{n},\text{ \qquad }x\in I_{n} \\ 
0,\,\,\,\text{\qquad }x\in \overline{I}_{n}%
\end{array}%
\right. ,  \label{dir1}
\end{equation}

Let $n\in\mathbb{N}$ and $n=\sum\limits_{i=0}^{\infty }n_{i}2^{i}.$ Then 
\begin{equation}
D_{n}\left( x\right) =w_{n}\left( x\right) \sum\limits_{j=0}^{\infty
}n_{j}w_{2^{j}}\left( x\right) D_{2^{j}}\left( x\right) .  \label{dir2}
\end{equation}

Set $G^2:= G\times G .$ The rectangular partial sums $ S_{M,N} $ of the 2-dimensional
Walsh-Fourier series of a function $f\in L_{2}\left( G^2\right) $ are defined
as follows:

\begin{equation*}
S_{M,N}f\left( x,y\right) :=\sum\limits_{i=0}^{M-1}\sum\limits_{j=0}^{N-1}%
\widehat{f}\left( i,j\right) w_{i}\left( x\right) w_{j}\left( y\right) ,
\end{equation*}
where the number 
\begin{equation*}
\widehat{f}\left( i,j\right) =\int\limits_{G^2}f\left( x,y\right)
w_{i}\left( x\right) w_{j}\left( y\right) d\mu \left( x,y\right)
\end{equation*}
is said to be the $\left( i,j\right) -$th Walsh-Fourier coefficient of the
function \thinspace $f.$

The norms (or quasi-norms) of the spaces $L_{p}(G^2)$ and $weak-L_{p}\left(
G^2\right) $ are respectively defined by 
\begin{equation*}
\left\| f\right\| _{p}:=\left( \int_{G^2}\left| f\right| ^{p}d\mu \right)
^{1/p},\qquad \left\| f\right\| _{weak-L_{p}(G^2)}:=\underset{\lambda >0}{%
\sup } \lambda \mu \left( f>\lambda \right) ^{1/p} \qquad \left( 0<p<\infty
\right) .
\end{equation*}

The $\sigma -$algebra generated by the dyadic 2-dimensional $I_{n}\left(
x\right) \times I_{n}\left( y\right) $ square of measure $2^{-n}\times
2^{-n} $ will be denoted by $\digamma _{n}\left( n\in \mathbb{N}\right) .$
Denote by $f=\left( f_{n},\text{ }n\in \mathbb{N}\right) $ the one-parameter
martingale with respect to $\digamma _{n}\left( n\in \mathbb{N}\right),$
(for details see e.g. Weisz \cite{Webook1} and \cite{Webook2}).
The maximal function of a martingale $f$ is defined by 
\begin{equation*}
f^{\ast }=\sup_{n\in \mathbb{N}}\left\vert f_{n}\right\vert .
\end{equation*}
Let $f\in L_{1}(G^2) $. Then the dyadic maximal function is given by 
\begin{equation*}
f^{*}\left( x,y\right) =\sup\limits_{n\in \mathbb{N}}\frac{1}{\mu \left(
I_{n}(x)\times I_{n}(y)\right) }\left| \int\limits_{I_{n}(x)\times
I_{n}(y)}f\left( s,t\right) d\mu \left( s,t\right) \right| ,\,\, \left(
x,y\right) \in G^2.
\end{equation*}
The dyadic Hardy martingale space $H_{p}(G^2)$ $\left( 0<p<\infty \right) $
consists of all functions $ f $ for which

\begin{equation*}
\left\| f\right\| _{H_{p}(G^2)}:=\left\| f^{*}\right\| _{p}<\infty .
\end{equation*}

If $f\in L_{1}(G^2) ,$ then it is easy to show that the sequence $\left(
S_{2^{n},2^{n}}f :n\in \mathbb{N}\right) $ is a martingale. If $f=\left(
f_{n},n\in \mathbb{N}\right) $ is a martingale, then the Walsh-Fourier
coefficients must be defined in a slightly different manner: $\qquad \qquad $
\begin{equation*}
\widehat{f}\left( i,j\right) :=\lim_{k\rightarrow \infty
}\int_{G}f_{k}\left( x,y\right) w_{i}\left( x\right) w_{j}\left( y\right)
d\mu \left( x,y\right) .
\end{equation*}

The Walsh-Fourier coefficients of $f\in L_{1}\left( G^2\right) $ are the
same as those of the martingale $\left( S_{2^{n},2^{n}}f:n\in \mathbb{N}%
\right) $ obtained from $f.$

A bounded measurable function $a$ is a p-atom, if there exists a
dyadic\thinspace 2-dimensional cube $I\times I\mathbf{,}$ such that 
\begin{equation*}
\int_{I\times I}ad\mu =0,\text{ \ \ }\left\Vert a\right\Vert _{\infty }\leq
\mu (I\times I)^{-1/p},\text{ \ \ \ supp}\left( a\right) \subset I\times I.
\end{equation*}

The dyadic Hardy martingale spaces $H_{p}\left( G^2\right) $ for $0<p\leq 1$ have atomic characterizations (for details see e.g. Weisz \cite{Webook1} and \cite{Webook2}).

\begin{lemma}
A martingale $f=\left( f_{n}:n\in \mathbb{N}\right) $ is in $H_{p}\left(
G^2\right) \left( 0<p\leq 1\right) $ if and only if there exist a sequence $%
\left( a_{k},k\in \mathbb{N}\right) $ of p-atoms and a sequence $\left(\mu
_{k},k\in \mathbb{N}\right) $ of a real numbers such that 
\begin{equation}
\qquad \sum_{k=0}^{\infty }\mu _{k}S_{2^n,2^n}a_{k}=f_{n}  \label{a1}
\end{equation}
and 
\begin{equation*}
\qquad \sum_{k=0}^{\infty }\left| \mu _{k}\right| ^{p}<\infty.
\end{equation*}
Moreover,   
$$\left\|f\right\|_{H_{p}}\backsim \inf\left( \sum_{k=0}^{\infty}\left|\mu_{k}\right|^{p}\right)^{1/p},$$ 
where the infimum is taken over
all decomposition of $f$ of the form (\ref{a1}).
\end{lemma}

It is known $\left[ 7,\text{ p. 125}\right] $ that the Walsh-Paley system is
not a Schauder basis in $L_{1}\left( G^2\right) $. Moreover, (for details
see \cite{S-W-S}) there exists a function in the dyadic martingale Hardy
space $H_{p}\left( G^2\right) $ $(0<p\leq 1)$ for which the partial sums are not
bounded in $L_{p}\left( G^2\right) .$ However, Weisz \cite{We} proved
the following estimation:

\textbf{Theorem W.} Let $\alpha \geq 0$, $0<p\leq1$ and $f\in H_{p}\left(
G^2\right).$ Then there exists an absolute constant $c_{p},$ depending only
on $p,$ such that

\begin{equation*}
\underset{n,m\geq 2}{\sup }\left( \frac{1}{\log n\log m}\right) ^{\left[ p%
\right] }\underset{2^{-\alpha }\leq k/l\leq 2^{\alpha },\text{ }\left(
k,l\right) \leq \left( n,m\right) }{\sum }\frac{\left\Vert
S_{k,l}f\right\Vert _{p}^{p}}{\left( kl\right) ^{2-p}}\leq c_{p}\left\Vert
f\right\Vert _{H_{p}}^{p},
\end{equation*}%
where $\left[ p\right] $ denotes the integer part of $p.$

In the case when $\alpha =0$ and $p=1$ from Theorem W we obtain that the
following is true:

\textbf{Theorem W1.} Let $f\in H_{1}\left( G^{2}\right) .$ Then there exists
an absolute constant $c,$ such that 
\begin{equation}
\frac{1}{\log ^{2}n}\underset{k=0}{\sum^{n}}\frac{\left\Vert
S_{k,k}f\right\Vert _{1}}{k^{2}}\leq c\left\Vert f\right\Vert _{H_{1}\left(
G^{2}\right) }.  \label{1}
\end{equation}

When $\alpha =0$ and $0<p<1$ from Theorem W follows the following result:

\textbf{Theorem W2.} Let $0<p<1$ and $f\in H_{p}\left( G^{2}\right) .$ Then
there exists an absolute constant $c_{p},$ depending only on $p,$ such that

\begin{equation}
\underset{k=0}{\sum^{\infty }}\frac{\left\Vert S_{k,k}f\right\Vert _{p}^{p}}{%
k^{4-2p}}\leq c_{p}\left\Vert f\right\Vert _{H_{p}\left( G^{2}\right) }^{p}.
\label{1.1}
\end{equation}

Goginava and Gogoladze \cite{gg} generalized the estimate (\ref{1}) (for details see \cite{MST} and \cite{tep7}) and proved that for any $f\in H_{1}\left( G^{2}\right) $, there exists an absolute constant $c$, such that 
\begin{equation}
\sum\limits_{n=1}^{\infty }\frac{\left\Vert S_{n,n}f\right\Vert _{1}}{n\log^{2}n}\leq c\left\Vert f\right\Vert_{H_{1}\left( G^{2}\right) }.
\label{th}
\end{equation}

The estimate (\ref{1.1}) was generalized in \cite{tep7} and it was proved that for any 
$0<p<1$ and $f\in H_{p}\left( G^{2}\right) $, there exists an absolute constant 
$c_{p}$, depending only on $p,$ such that 
\begin{equation}
\sum\limits_{n=1}^{\infty }\frac{\left\Vert S_{n,n}f\right\Vert _{p}^{p}}{%
n^{3-2p}}\leq c_{p}\left\Vert f\right\Vert _{H_{p}\left( G^{2}\right) }^{p}.
\label{th1}
\end{equation}

In \cite{tep6} and \cite{tep7} it was proved that the following is true:

\textbf{Theorem T1.} Let $0<p\leq 1$ and $\Phi :\mathbb{N}\rightarrow
\lbrack 1,$ $\infty )$ be any nondecreasing function, satisfying the
condition $\underset{n\rightarrow \infty }{\lim }\Phi \left( n\right)
=+\infty .$ Then there exists a martingale $f\in H_{p}\left( G^{2}\right) $
such that

\begin{equation*}
\underset{n=1}{\overset{\infty }{\sum }}\frac{\left\Vert S_{n,n}f\right\Vert
_{weak-L_{p}}^{p}\Phi \left( n\right) }{n^{3-2p}\log ^{2[p]}n}=\infty ,
\end{equation*}%
where $\left[ p\right] $ denotes the integer part of $p.$

Analogical problems for the one-dimensional case were proved in \cite{bnpt1}, \cite{gat1}, \cite{si1}, \cite{tep8}, \cite{tep9}, \cite{tut1} and for the two-dimensional case in \cite{MST}, \cite{n1}, \cite{PTW3}.

In view of Theorem T1 we can conclude that the sequence $\left\{ 1/(n^{3-2p}\log
^{2[p]}n)\right\} _{n=1}^{\infty }$ in inequalities (\ref{th}) and (\ref{th1}%
) can not be improved, which gives sharpness for $\alpha =0$.

In this paper we consider the analogous problem for $\alpha>0$ and prove the sharpness of the sequence in Theorem W when $\alpha >0$:

\begin{theorem}
Let $0<p<1,$ $\alpha >0$ and $\Phi :\mathbb{N}^{2}\rightarrow \lbrack 1,$ $%
\infty )$ satisfies the conditions
\begin{equation}
\Phi \left( m_{1},n_{1}\right) \geq \Phi \left( m_{2},n_{2}\right) \text{
when }m_{1}\geq m_{2}\text{ \ and }n_{1}\geq n_{2}  \label{dir111}
\end{equation}%
and 
\begin{equation}
\underset{m,n\rightarrow \infty }{\lim }\Phi \left( m,n\right) =+\infty .
\label{dir222}
\end{equation}%
Then there exists a martingale $f\in H_{p}\left( G^{2}\right) $ such that

\begin{equation*}
\underset{n,m\in \mathbb{N}_+}{\sup}\underset{2^{-\alpha}\leq k/l\leq
2^{\alpha },\text{ }\left(k,l\right)\leq\left(n,m\right)}{\sum}\frac{%
\left\Vert S_{k,l}f\right\Vert_{p}^{p}\Phi\left(k,l\right)}{%
\left(kl\right)^{2-p}} =\infty.
\end{equation*}
\end{theorem}

\textbf{Proof.} Under the condition (\ref{dir222}) there exists an
increasing sequence of positive integers $\left\{ \alpha _{k}:\text{ }%
k\geq 0\right\} $ such that $\alpha _{0}\geq 2,$ \ $\alpha_k+[\alpha]+1<\alpha_{k+1}$ \ and 
\begin{equation}
\sum_{k=0}^{\infty }{\Phi ^{-p/4}\left( 2^{\alpha _{k}},2^{\alpha
_{k}}\right) }<\infty .  \label{2}
\end{equation}

Let \qquad 
\begin{equation*}
f_{n}=\sum_{\left\{ k;\text{ }\alpha _{k}+[\alpha]+1<n\right\} }\lambda
_{k}a_{k},
\end{equation*}
where 

\begin{equation*}
\lambda _{k}:=2^{2[\alpha]+2}{\Phi ^{-1/4}\left( 2^{\alpha _{k}},2^{\alpha _{k}}\right) }
\end{equation*}
and
\begin{equation*}
a_{k}\left(x,y\right):=2^{\alpha _{k}\left( 2/p-2\right)-2[\alpha]-2}\left(
D_{2^{\alpha_k+[\alpha]+1}}\left(x\right) -D_{2^{\alpha _{k}}}\left(
x\right) \right) \left( D_{2^{\alpha _{k}+[\alpha]+1}}\left( y\right)
-D_{2^{\alpha _{k}}}\left( y\right) \right) .
\end{equation*}

Since
\begin{equation*}
S_{2^{n},2^{n}}a_{k}=\left\{ 
\begin{array}{l}
a_{k}\text{, \qquad } \alpha _{k}+[\alpha]+1<n, \\ 
0\text{, \qquad } \alpha _{k}+[\alpha]+1\geq n,%
\end{array}
\right.
\end{equation*}

\begin{eqnarray*}
&&\text{supp}(a_{k})=I_{\alpha _{k}}^2, \text{ \qquad } \ \ \ \int_{I_{\alpha
_{k}}^2}a_{k}d\mu=0, \text{ \qquad } \ \ \ \left\| a_{k}\right\| _{\infty }\leq
2^{2\alpha _{k}/p}=(\mu(\text{supp}a_{k}))^{-1/p}
\end{eqnarray*}
from Lemma 1 and (\ref{2}) we obtain that $f\in H_{p}(G^2).$

It is obvious that 

\begin{eqnarray}  \label{5}
&&\widehat{f}(i,j) \\
&=&\left\{ 
\begin{array}{l}
\frac{2^{\alpha _{k}\left( 2/p-2\right) }}{\Phi ^{1/4}\left( 2^{\alpha
_{k}},2^{\alpha _{k}}\right)}, \text{ \ if \ }\left( i,j\right) \in
\left\{2^{\alpha _{k}},...,2^{\alpha _{k}+[\alpha]+1}-1\right\} \times
\left\{ 2^{\alpha _{k}},...,2^{\alpha _{k}+[\alpha]+1}-1\right\} ,\text{ }%
k\in \mathbb{N} \\ 
0,\text{ \ if \ }\left(i,j\right) \notin \bigcup\limits_{k=1}^{\infty
}\left\{ 2^{\alpha _{k}},...,2^{\alpha _{k}+[\alpha]+1}-1\right\}\times
\left\{ 2^{\alpha _{k}},...,2^{\alpha _{k}+[\alpha]+1}-1\right\} .%
\end{array}%
\right.  \notag
\end{eqnarray}

Let\textbf{\ } $2^{\alpha _{k}}<m,n<2^{\alpha _{k}+[\alpha ]+1}$. In view of (\ref{5}) we can conclude that 

\begin{eqnarray}
&&S_{m,n}f\left( x,y\right)  \label{13d} \\
&=&\sum_{\eta =0}^{k-1}\sum_{i=2^{\alpha _{\eta }}}^{2^{\alpha _{\eta
}+1}-1}\sum_{j=2^{\alpha _{\eta }}}^{2^{\alpha _{\eta }+1}-1}\widehat{f}%
(i,j)w_{i}\left( x\right) w_{j}\left( y\right)+\sum_{i=2^{\alpha _{k}}}^{m-1}\sum_{j=2^{\alpha _{k}}}^{n-1}\widehat{f}%
(i,j)w_{i}\left( x\right) w_{j}\left( y\right)  \notag \\
&=&\sum_{\eta =0}^{k-1}\sum_{i=2^{\alpha _{\eta }}}^{2^{\alpha _{\eta
}+1}-1}\sum_{j=2^{\alpha _{\eta }}}^{2^{\alpha _{\eta }+1}-1}\frac{2^{\alpha
_{\eta }\left( 2/p-2\right) }w_{i}\left( x\right) w_{j}\left( y\right) }{%
\Phi ^{1/4}\left( 2^{\alpha _{\eta }},2^{\alpha _{\eta }}\right) }  +\sum_{i=2^{\alpha _{k}}}^{m-1}\sum_{j=2^{\alpha _{k}}}^{n-1}\frac{%
2^{\alpha _{k}\left( 2/p-2\right) }w_{i}\left( x\right) w_{j}\left( y\right) 
}{\Phi ^{1/4}\left( 2^{\alpha _{k}},2^{\alpha _{k}}\right) }  \notag \\
&=&\sum_{\eta =0}^{k-1}\frac{2^{\alpha _{\eta }\left( 2/p-2\right) }}{\Phi
^{1/4}\left( 2^{\alpha _{\eta }},2^{\alpha _{\eta }}\right) }\left(
D_{2^{\alpha _{\eta }+1}}\left( x\right) -D_{2^{\alpha _{\eta }}}\left(
x\right) \right) \left( D_{2^{\alpha _{\eta }+1}}\left( y\right)
-D_{2^{\alpha _{\eta }}}\left( y\right) \right)  \notag \\
&&+\frac{2^{\alpha _{k}\left( 2/p-2\right) }}{\Phi ^{1/4}\left( 2^{\alpha
_{k}},2^{\alpha _{k}}\right) }\left( D_{m}\left( x\right) -D_{2^{\alpha
_{k}}}\left( x\right) \right) \left( D_{n}\left( y\right) -D_{2^{\alpha
_{k}}}\left( y\right) \right)  \notag \\
&:=&I+II.  \notag
\end{eqnarray}

Let $\left( x,y\right) \in \left( G\backslash I_{1}\right) \times \left(
G\backslash I_{1}\right),$ $n$ and $m$ are odd numbers, such that $2^{\alpha _{k}}<m,n<2^{\alpha _{k}+[\alpha ]+1}$. Since  $\alpha _{k}\geq 2$ $\left( k\in \mathbb{N}\right),$ if we combine (\ref{dir1}) and (\ref{dir2}) it follows that
$$D_{2^{\alpha_k}}\left( x\right)=D_{2^{\alpha_k}}\left( y\right)=0$$
and
\begin{eqnarray}  \label{13a}
\left\vert II\right\vert=\frac{2^{\alpha_k\left( 2/p-2\right)}}{\Phi ^{1/4}\left(2^{\alpha_k},2^{\alpha_{k}}\right)}\left\vert w_m\left( x\right)w_1\left(x\right)D_1\left(x\right) w_n \left( y\right)D_1\left( y\right)\right\vert=\frac{2^{\alpha _{k}\left(2/p-2\right)}}{\Phi ^{1/4}\left(2^{\alpha_k},2^{\alpha_k}\right)}. 
\end{eqnarray}
By applying (\ref{dir1}) and the condition $\alpha _{n}\geq 2$ $\left( n\in 
\mathbb{N}\right) $ for $I$ we have that

\begin{equation}
I=\sum_{\eta =0}^{k-1}\frac{2^{\alpha _{\eta}\left( 2/p-2\right) }}{\Phi
^{1/4}\left( 2^{\alpha _{\eta }}\right) }\left( D_{2^{\alpha _{\eta
}+1}}\left( x\right) -D_{2^{\alpha _{\eta }}}\left( x\right) \right) \left(
D_{2^{\alpha _{\eta }+1}}\left(y\right)-D_{2^{\alpha _{\eta }}}\left(
y\right) \right) =0.  \label{13b}
\end{equation}
By combining (\ref{13a}) and (\ref{13b}), for $2^{\alpha _{k}}<m,n<2^{\alpha _{k}+[\alpha ]+1}$ we get that

\begin{eqnarray} \label{13}
&&\left\Vert S_{m,n}f \right\Vert _{weak-L_{p}}  
\\
&\geq &\frac{2^{\alpha _{k}\left( 2/p-2\right) }}{2\Phi ^{1/4}\left(
2^{\alpha _{k}},2^{\alpha _{k}}\right) }\left( \mu \left\{ \left( x,y\right)
\in \left( G\backslash I_{1}\right) \times \left( G\backslash I_{1}\right)
:\left\vert S_{m,n}f\left( x,y\right) \right\vert \geq \frac{2^{\alpha
_{k}\left( 2/p-2\right) }}{2\Phi ^{1/4}\left( 2^{\alpha _{k}},2^{\alpha
_{k}}\right) }\right\} \right) ^{1/p}  \notag \\
&\geq &\frac{2^{\alpha _{k}\left( 2/p-2\right) }}{2\Phi ^{1/4}\left(
2^{\alpha _{k}},2^{\alpha _{k}}\right) }\mu \left( \left( G\backslash
I_{1}\right) \times \left( G\backslash I_{1}\right) \right) \geq \frac{%
c_{p}2^{\alpha _{k}\left( 2/p-2\right) }}{\Phi ^{1/4}\left( 2^{\alpha
_{k}},2^{\alpha _{k}}\right) }.  \notag
\end{eqnarray}

According to (\ref{dir111}), (\ref{dir222}) and (\ref{13}) we can conclude
that

\begin{eqnarray} \label{14}
&&\underset{s,l\in \mathbb{N}_{+}}{\sup }\underset{2^{-\alpha }\leq m/n\leq
2^{\alpha },\text{ }\left( m,n\right) \leq \left( s,l\right) }{\sum }\frac{%
\left\Vert S_{m,n}f\right\Vert _{weak-L_{p}}^{p}\Phi \left( m,n\right) }{%
\left( mn\right) ^{2-p}} \\ \notag
&\geq &\underset{2^{\alpha _{k}}<m,n\leq 2^{\alpha _{k}+\alpha }}{\sum }%
\frac{\left\Vert S_{m,n}f\right\Vert _{weak-L_{p}}^{p}\Phi \left( m,n\right) 
}{\left( mn\right) ^{2-p}} \\ \notag
&\geq &\frac{c_{p}\Phi \left( 2^{\alpha _{k}},2^{\alpha _{k}}\right) }{%
2^{\alpha _{k}\left( 4-2p\right) }}\underset{2^{\alpha _{k}}<m,n\leq
2^{\alpha _{k}+\alpha }}{\sum }\left\Vert S_{m,n}f\right\Vert
_{weak-L_{p}}^{p} \\ \notag
&\geq &\frac{c_{p}\Phi \left( 2^{\alpha _{k}},2^{\alpha _{k}}\right) }{%
2^{\alpha _{k}\left( 4-2p\right) }}\underset{2^{\alpha _{k}-1}<m,n\leq
2^{\alpha _{k}-1+\alpha }}{\sum }\left\Vert S_{2m+1,2n+1}f\right\Vert
_{weak-L_{p}}^{p} 
\end{eqnarray}
\begin{eqnarray*}
&\geq &\frac{c_{p}\Phi \left( 2^{\alpha _{k}},2^{\alpha _{k}}\right) }{%
2^{\alpha _{k}\left( 4-2p\right) }}\frac{2^{\alpha _{k}\left(2-2p\right) }}{%
\Phi ^{1/4}\left( 2^{\alpha _{k}},2^{\alpha _{k}}\right) }\underset{%
2^{\alpha _{k}-1}<m,n\leq 2^{\alpha _{k}-1+\alpha }}{\sum }1 \\
&\geq &\frac{c_{p}\Phi ^{3/4}\left( 2^{\alpha _{k}},2^{\alpha _{k}}\right) }{%
2^{2\alpha _{k}}}\left( 2^{\alpha _{k}-1+\alpha }-2^{\alpha _{k}-1}-1\right)
^{2} \\
&\geq &c_{p}\Phi ^{3/4}\left( 2^{\alpha _{k}},2^{\alpha _{k}}\right)
\rightarrow \infty,\qquad\text{as}\text{\qquad }k\rightarrow \infty .
\end{eqnarray*}

By combining (\ref{2}-\ref{14}) we complete the proof of Theorem 1.

\end{document}